\documentclass[12pt]{article}
\usepackage{fullpage}
\usepackage{amstex}
\usepackage{amssymb}
\usepackage{theorem}

{\theorembodyfont{\rmfamily}
}

\def\MOD{\mathop{\rm mod}}
\def\tet{{\theta}}
\def\l{{\langle}}
\def\r{{\rangle}}
\def\del{\delta}
\def\ome{\omega}

\def\bks{\backslash}	
\def\tilA{\widetilde A}	
\def\tilc{\widetilde c}	
\def\tilj{\widetilde j}	
\def\tilrho{\widetilde\rho}	
\def\tilM{\widetilde M}	
\def\tilnu{\widetilde\nu}	
\def\onu{\overline\nu}	
\def\hb{\hfill\break}
\def\gam{\gamma}	
\def\kap{\kappa}	
\def\lam{\lambda}

\def\ome{\omega}
\def\nek{,\ldots,}
\def\sig{\sigma}
\def\alp{\alpha}
\def\bet{\beta}

\def\prodl{\prod\limits}

\title{On Gaps under GCH Type Assumptions}
\author{Moti Gitik\\
School of Mathematical Sciences\\
Tel Aviv University\\
 Tel Aviv 69978, Israel\\
gitik@@ math.tau.ac.il}
\date{}
\begin{document}
\maketitle
\begin{abstract}
We prove equiconsistency results concerning gaps between
a singular strong limit cardinal $\kappa$ of
cofinality $\aleph_0$ and its power under
assumptions that $2^\kap =\kap^{+\del +1}$  for 
$\del <\kap$  and some weak form of the
Singular Cardinal Hypothesis below $\kap$.
Together with previous results this basically
completes the study of consistency strength of 
various gaps between such $\kap$  and its power
under GCH type assumptions below.
\end{abstract}

\baselineskip=18pt
\setcounter{section}{-1}
\renewcommand{\labelenumi}{(\theenumi)}
\renewcommand{\labelenumiii}{(\theenumiii)}
\section{Introduction}

Our first result deals with cardinal gaps.

We continue [Git-Mit] and show the
following:

\medskip
\noindent
{\bf Theorem 1.}\quad  {\sl Suppose that $\kap$  is a
strong limit cardinal of cofinality $\aleph_0$,
$\del <\kap$ is a cardinal of uncountable cofinality.
If $2^\kap\ge\kap^{+\del}$  and the
Singular Cardinal Hypothesis holds below $\kap$ at least
for cardinals of cofinality $cf\del$, then in the core
model either
\begin{itemize}
\item [(i)] $o(\kap)\ge\kap^{+\del+1}+1$ or
\smallskip
\item [(ii)] $\{\alp <\kap\mid o(\alp)\ge
\alp^{+\del+1}+1\}$ is unbounded in $\kap$.
\end{itemize}
}

Together with [Git-Mag] and [Git1] this provides
the equiconsistency result for cardinal gaps of
uncountable cofinality. 
Surprisingly the proof
uses very little of the indiscernibles theory for
extenders developed in [Git-Mit].  Instead, basic
results of the Shelah pcf-theory play the
crucial role.

Building on the analysis of indiscernibles for
uncountable cofinality of [Git-Mit] and
pcf-theory we show the following:

\medskip
\noindent
{\bf Theorem 2.}\quad {\sl If for a set $a$  of
regular cardinals above $2^{|a|^++\aleph_2}$
$|pcf a|>|a|+\aleph_1$ then there is an inner
model with a strong cardinal.}  

Using this result, we extend Theorem 1 to ordinal gaps:

\medskip
\noindent
{\bf Theorem 3.}\quad  {\sl Suppose that $\kap$  is a
strong limit cardinal of cofinality $\aleph_0$,
$\del <\kap$ is a cardinal above $\aleph_1$ of
uncountable cofinality and $\ell <\ome$.  If
$2^\kap\ge\kap^{+\del^\ell}$  and the
Singular Cardinal Hypothesis holds below $\kap$ at least
for cardinals of cofinality $cf\del$, then in the core
model either
\begin{itemize}
\item [(i)] $o(\kap)\ge\kap^{+\del^\ell+1}+1$ or
\smallskip
\item [(ii)] $\{\alp <\kap\mid o(\alp)\ge
\alp^{+\del^\ell+1}+1\}$ is unbounded in $\kap$.
\end{itemize}
}

If the pcf structure between $\kap$  and $2^\kap$
is not ``wild" (thus, for example, if there is
no measurable of the core model between $\kap$
and $2^\kap$), then the result holds also for
$\del=\aleph_1$.

These theorems and related results are proved in
Section 1 of the paper.
Actually more general results (1.20, 1.21) are
proved for ordinal gaps but the formulations
require technical notions ``Kinds" and
``Kinds$^*$" and we will not reproduce them
here.
In Section 2 we sketch some complimentary
forcing constructions based on [Git1].  Thus we
are able to deal with cardinal gaps of
cofinality $\aleph_0$  and show the following
which together with Theorem 1 provides the
equiconsistency for the cases of cofinality $\aleph_0$.

\medskip
\noindent
{\bf Theorem 4.}\quad {\sl Suppose that in the core
model $\kap$ is a singular cardinal of
cofinality $\aleph_0$,  $\del <\kap$  is a
cardinal of cofinality $\aleph_0$  as well and
for every $\tau <\del$ the set $\{\alp <\kap\mid
o(\alp)\ge\alp^{+\tau}\}$  is unbounded in
$\kap$.  {\it Then\/} for every $\alp <\del^+$
there is a cofinalities preserving, not adding
new bounded subsets to $\kap$  extension satisfying
$2^\kap\ge\kap^{+\alp}$.}

The Rado-Milner paradox is used to show the
following:

\medskip
\noindent
{\bf Theorem 5.}\quad {\sl Suppose that in the core
model $\kap$ is a singular cardinal of
cofinality $\aleph_0$, $\del <\kap$  is a cardinal
of uncountable cofinality and for every $n<\ome$
the set $\{\alp <\kap\mid o(\alp)\ge\alp^{+\del^n}\}$
is unbounded in $\kap$.  {\it Then\/} for every
$\alp <\del^+$  there is cofinalities preserving
not adding new bounded subsets to $\kap$
extension satisfying $2^\kap\ge\kap^{+\alp}$.}

A more general result (2.6) of the same flavor is
obtained for ordinal gaps.

In the last section, we summarize the situation and
discuss related open questions and some further
directions.

A knowledge of the basic $pcf$-theory results is
needed for Section 1.  We refer to the Burke-Magidor
[Bur-Mag] survey paper or to Shelah's book
[Sh-g] on these matters.  Results on ordinal
gaps and the strength of ``$|pcf a|>|a|$"
require in addition familiarity with basics of
indiscernible structure for extenders.  See
Gitik-Mitchell [Git-Mit] on this subject. 

Results of Sections 2 are built on short
extender based Prikry forcings, mainly those of
[Git1]. 

\medskip
\noindent
{\bf Acknowledgement.}\quad We are grateful to
Saharon Shelah for many helpful conversations
and for explanations that he gave on the pcf-theory.

\section{On the Strength of Gaps}

Let $SSH^\del_{<\kap}$  ($SSH^{\le\del}_{<\kap}$)
denote the Shelah Strong
Hypothesis below $\kap$ for cofinality $\del$ $(\le\del)$
which means that for every singular cardinal
$\tau <\kap$  of cofinality $\del(\le\del)$ $\ pp(\tau)=\tau^+$.
We assume that there is no inner model with a strong
cardinal.  First we will prove the following:

\medskip
\noindent
{\bf Theorem 1.1.}\quad {\sl Suppose that $\kap$
is a singular strong limit cardinal of
cofinality $\aleph_0$, $\del <\kap$  a cardinal
of uncountable cofinality, $2^\kap\ge\kap^{+\del}$  and
$SSH^{cf\del}_{<\kap}$.  Then in the core model either 
\begin{enumerate}
\item [(i)] $o(\kap)\ge\kap^{+\del +1}+1$
\smallskip\noindent
or
\smallskip
\item [(ii)] $\{\alp <\kap\mid o(\alp)\ge\alp^{+\del +1}
+1\}$
is unbounded in $\kap$.
\end{enumerate}
}

\medskip
\noindent
{\bf Remark 1.2.}\quad (1) in either case we have
in the core model a cardinal $\alp$  carrying an 
extender of the length $\alp^{+\del +1}$.

(2) By [Git-Mag] or [Git1] it is  possible to
force, using (i) or (ii), the situation assumed
in the theorem.  So this provides equiconsistency
result.

\medskip
\noindent
{\bf Proof.}\quad If $\del$  is a regular cardinal then
let $A$  be the set of cardinals $\kap^{+\tau +1}$  so
that $\tau <\del$  and either $o(\alp)<\kap^{+\tau}$  for
every $\alp <\kap^{+\tau}$ or else $\kap^{+\tau}$  is
above every measurable of the core model smaller than
$\kap^{+\del}$.  The set $A$  is unbounded in
$\kap^{+\del}$  since there is no overlapping
extenders in the core model.  If $cf\del <\del$  then we
fix $\l\del_i\mid i<cf\del\r$ an increasing
sequence of regular cardinals with limit $\del$.  
For every $i<cf\del$  define $A_i$  to be the
set of cardinals $\kap^{+\tau +1}$  so that
$\tau <\del_i$  and either
$o(\alp)<\kap^{+\tau}$  for every $\alp
<\kap^{+\tau}$  or else $\kap^{+\tau}$  is above
every measurable of the core model smaller than
$\kap^{+\del_i}$.  Again, each of $A_i$'s will
be unbounded in $\kap^{+\del_i}$ since there is
no overlapping extenders in the core model.

The following fact was proved in [Git-Mit,
3.24]:	

\medskip
\noindent
{\bf Claim 1.3.}\quad If $B\subseteq A$  in case
$cf\del =\del$ or $B\subseteq A_i$  for some
$i<cf\del$, in case $cf\del <\del$ then $|B|<\inf B$
implies $\max(pcf(B))=(\sup B)^+$.

Now for every $\kap^{+\alp +1}\in A$  or $\kap^{+\alp
+1}\in\bigcup_{i<cf\del}A_i$  (if $cf\del <\del$) we pick
a set $\{c^\alp_n\mid n<\ome\}$  of regular cardinals
below $\kap$  so that $\kap^{+\alp+1}\in pcf\{c^\alp_n
\mid n<\ome\}$.
Set
$$a=\{c^\alp_n\mid n<\ome\ ,\quad\kap^{+\alp +1}\in A
\quad\hbox{if}\quad cf\del =\del\quad\hbox{or}\quad
\kap^{+\alp+1}\in\bigcup_{i <cf\del}A_i\quad
\hbox{otherwise}\}\ .$$
Removing its bounded part, if necessary, we can
assume that $\min a>|a|^+$. 

\medskip
\noindent
{\bf Claim 1.4.}\quad For every $b\subseteq a\mid
A\cap pcf (b)|\le |b|$ or $|A_i\cap pcf (b)|\le
|b|$, for every $i<cf\del$, if $cf\del <\del$.

\medskip
\noindent
{\bf Proof.}\quad It follows from Shelah's Localization
Theorem [Sh-g] and Claim 1.3.\hfill $\square$

In particular, $|a|=\del$.

Let $b_{\kap^+}[a]$  be the pcf-generator corresponding
to $\kap^+$.  Consider $a^*=a\bks b_{\kap^+}[a]$.
For every $\alp >0$,  if $\kap^{+\alp +1}\in A$
or $\bigcup_{i<cf\del}A_i$ then $\kap^{+\alp +1}\in
pcf(a^*)$.  Hence, $|(pcf a^*)\cap A|=\del$  or
$|pcf(a^*)\cap A_i\mid =\del_i$ for each
$i<cf\del$ and by Claim 1.4, then $|a^*|=\del$.

\medskip
\noindent
{\bf Claim 1.5.}\quad Let $\l\tau_n\mid
n<\ome\r$ be an increasing unbounded in $\kap$
sequence of limit points of $a^*$ of cofinality
$cf\del$.  Then for every ultrafilter $D$  on
$\ome$ including all cofinite sets\hb
$cf\Big(\prodl_{n<\ome}\tau^+_n/D\Big)>\kap^+$.

\medskip
\noindent
{\bf Proof.}\quad For every $n<\ome$, $\tau_n$
is a singular cardinal of cofinality $cf\del$.
So, by the assumption $pp(\tau_n)=\tau^+_n$.
Then $\tau^+_n=cf (\prod t/E)$, for every
unbounded in $\tau_n$  set of regular cardinals
with $|t|<\tau_n$  and an ultrafilter $E$
on it including all cobounded subset of $t$. In
particular, $\tau^+_n\in pcf (a^*\cap\tau_n)$
since $\tau_n$  is a limit point of $a^*$. 

So $\{\tau^+_n\mid n<\ome\}\subseteq pcf a^*$.
By [Sh-g], then $pcf \{\tau^+_n\mid n<\ome\}\subseteq
pcf (pcf a^*)=pcf a^*$.  But by the choice of
$a^*$,  $\kap^+\notin pcf a^*$.  Hence for every
ultrafilter $D$  on $\ome$, cf
$\Big(\prodl_{n<\ome}\tau^+_n/D\Big)\not= \kap^+$.
\hfill $\square$

Now, $|a^*|=\del,\cup a^*=\kap,\ cf\del
>\aleph_0$ and $cf\kap =\aleph_0$.  Hence there
is an increasing unbounded in $\kap$  sequence
$\l\tau_n\mid n <\ome\r$ of limit
points of $a^*$  so that for every $n>0$
$|a^*\cap (\tau_{n-1},\tau_n)|=\del$  and
$|(a^*\cap\tau_n)\bks \bet|=\del$  for every
$\bet <\tau_n$.  By Claim 1.5, $\l\tau^+_n\mid
n<\ome\r$  are limits of indiscernibles.
We refer to [Git-Mit] for basic facts on this
matter used here.  There is a principal
indiscernible $\rho_n\le\tau^+_n$ for all but
finitely many $n$'s.  By the Mitchell Weak Covering Lemma,
$\tau^+_n$  in the sense of the core model is
the real $\tau^+_n$,  since $\tau_n$  is
singular.  This implies that $\rho_n\le\tau_n$,
since a principal indiscernible cannot be
successor cardinal of the core model.  Also,
$\rho_n$  cannot be $\tau_n$,  since again
$\tau^+_n$ computed in the core model correctly
and so there is no indiscernibles between
measurable now $\tau_n$  and its successor
$\tau^+_n$.  Hence $\rho_n<\tau_n$.  By the
choice of $\tau_n$, the interval $(\rho_n,\tau_n)$
contains at least $\del$  regular cardinals.  So
$\rho_n$ is a principal indiscernible of
extender including at least $\del +1$  regular
cardinals which either seats over $\kap$  or
below $\kap$.  This implies that either
$o(\kap)\ge\kap^{+\del +1}+1$  or $\{\alp
<\kap\mid o(\alp)\ge\alp^{+\del +1}+1\}$  is
unbounded in $\kap$.\hfill$\square$  

Using the same ideas, let us show the following
somewhat more technical result:

\medskip
\noindent
{\bf Theorem 1.6.}\quad {\sl Let $\kap
=\bigcup_{n<\ome}\kap_n$  be a strong limit
cardinal with $\kap_0<\kap_1<\cdots
<\kap_n<\cdots$.

Assume $2^\kap\ge\kap^{++}$ and $SSH^{\aleph_1}_{<\kap}$  
(Shelah Strong Hypothesis below $\kap$ for
cofinality $\aleph_1$, i.e. $pp\tau =\tau^+$
for every singular $\tau <\kap$  of cofinality
$\aleph_1$).  Then there are at most countably many
principal indiscernibles $\l\rho_{n,m}\mid
m,n <\ome\r$  with indiscernibles $\l\del_{n,m}
\mid m,n<\ome\r$ so that for each $n,m<\ome$  
$\kap_n\le\rho_{n,m}\le\del_{n,m},\rho_{n,m}$ is
the principal indiscernible of $\del_{n,m}$,
each $\del_{n,m}$  is a regular cardinal and for
every $m<\ome$
$cf\Big(\prodl_{n<\ome}\del_{n,m}\Big/ D_m\Big)>\kap^+$,
where $D_m$  is an ultrafilter on $\ome$
including all cofinite sets.}

\medskip
\noindent
{\bf Remark 1.7.}\quad The theorem implies
results of the following type proved in
[Git-Mit]: if $2^\kap =\kap^{+m}$  $(2<m<\ome)$
and GCH below $\kap$,  then $o(\kap)\ge\kap^{+m}+1$,
provided that for some $k<\ome$  the set of
$\nu<\kap$  such that $o(\nu)>\nu^{+k}$ is
bounded in $\kap$. 

\medskip
\noindent
{\bf Proof.}\quad Suppose otherwise. 

Collapsing if necessary $2^\kap$  to
$\kap^{++}$, we can assume that $2^\kap
=\kap^{++}$.  Let $\l\rho_{n,i}\mid n
<\ome$, $i<\ome_1\r$  and
$\l\del_{n,i}\mid n<\ome, i<\ome_1\r$
witness the failure of the theorem.  We can
assume that for every $n<\ome$  and
$i<j<\ome_1$  
$$\rho_{n,i}\le\del_{n,i}<\rho_{n,j}\le \del_{nj}\ .$$
Let $a=\l\del_{n,i}\mid n<\ome, i<\ome_1\r$.
Consider $a^*=a\bks b_{\kap^+}[a]$.  Then for
every $i<\ome_1$ the set
$c_i=a^*\cap\{\del_{n,i}\mid n<\ome\}$ is
infinite, since
$cf\Big(\prodl_{n<\ome}\del_{n,i}\Big/D_i\Big)=\kap^{++}$
for some $D_i$.

The following is obvious. 

\medskip
\noindent
{\bf Claim 1.8.}\quad
There is an infinite set $d\subseteq\ome$  such that
for every $n\in d$  there are uncountably many
$i$'s with $\del_{n,i}\in c_i$.

For every $n\in d$ let
$$\tau_n=\sup\{\del_{n,i}\mid\del_{n,i}\in
C_i\}\ .$$
Then each such $\tau_n$ is a singular cardinal
of uncountable cofinality. 
Also, $\tau^+_n\in pcf a^*$  for every $n\in d$,
since $pp\tau_n=\tau^+_n$.  But then
$pcf\{\tau^+_n\mid n\in d\}\subseteq pcf a^*$.
Hence $\kap^+\not\in pcf \{\tau^+_n\mid n\in
d\}$.  Now, this implies as in the proof of 1.1
that $\tau^+_n$'s are indiscernibles and there
are principal indiscernibles for $\tau^+_n$'s
below $\tau_n$.  Here this is impossible since
then there should be overlapping extenders.
Contradiction.\hfill $\square$

We will use 1.6 further in order to deal with
ordinal gaps.

As above, we show the following
assuming that there is no inner model with a
strong cardinal.

\medskip
\noindent
{\bf Proposition 1.9.}\quad {\sl Suppose that $\l
\tau_\alp\mid \alp <\tet\r$  is an increasing
sequence of regular cardinals.  $\tet$ is a
regular cardinal  $>\aleph_1$  and
$\tau_0 >2^\tet$.  Then there is an unbounded
$S\subseteq\tet$  such that for every $\del$  of
uncountable cofinality which is a limit of
points of $S$ the following holds:

${\rm (*)}$ for every ultrafilter $D$  on $\del\cap S$
including all cobounded subsets of $\del\cap S$
$$tcf\Big(\prodl_{\alp\in\del\cap S}\tau_\alp\Big/
D\Big)=tcf\Big(\prodl_{\alp\in\del\cap S}\tau_\alp
\Big/J^{bd}_{\del\cap S}\Big)<\tau_{\alp +1}$$
where $J^{bd}_{\del\cap S}$ denotes the ideal of
bounded subsets of $\del\cap S$.
}

\medskip
\noindent
{\bf Proof.}\quad Here we apply the analysis of
indiscernibles of [Git-Mit] for uncountable
cofinality.  Let
$\l\nu_\bet\mid\bet\le\tet\r$ be the
increasing enumeration of the closure of
$\l\tau_\alp\mid\alp <\tet\r$.  Let
$A\subseteq\tet$  be the set of indexes of all
principal indiscernibles for $\nu_\tet$  among
$\nu_\bet$'s $(\bet <\tet)$.  Then $A$  is a
closed subset of $\tet$.  Now split into two
cases.

\medskip
\noindent
{\bf Case 1.}\quad $A$ is bounded in $\tet$.

Let $\bet^*=\sup A$.  We have a club
$C\subseteq\tet$  so that for every $\alp\in C$,
$\bet\in (\bet^*,\alp)$  if $\nu_\bet$  is a principal
indiscernible, then it is a principal
indiscernible for an ordinal below $\nu_\alp$.
Now let $\alp$  be a limit point of $C$  of uncountable
cofinality.  Then by results of [Git-Mit],
$pp\nu_\alp=\nu^+_\alp$ and moreover
$tcf\Big(\prodl_{\bet <\alp}\nu_\bet\Big/J^{bd}_{\nu_\alp}
\Big)=\nu^+_\alp$.  So we are
done.

\medskip
\noindent
{\bf Case 2.}\quad $A$ is bounded in $\tet$.

Let $\tilA$  be the set of limit points of $A$.
For every $\alp\in\tilA$  we consider 
$\nu_{\alp +1}$.  Let $\nu^*_{\alp +1}$  be the
principal indiscernible of $\nu_{\alp +1}$.
Then $\nu_\alp\le\nu^*_{\alp +1}\le\nu_{\alp +1}$.

The following is the main case:

\medskip
\noindent
{\bf Subcase 2.1.}\quad For every $\alp$  in an unbounded
set $S\subseteq\tet$, $\nu^*_{\alp +1}$  is a
principal indiscernible for $\nu_{\tet}$  and
$\nu_{\alp +1}$ is an indiscernible belonging to
some $\onu_{\alp +1}$  over $\nu_{\tet}$ of
cofinality $\ge\nu_{\ome_1}$ in the core model. 

We consider the set $B=\{\onu_{\alp +1}|\alp\in S\}$.
If $|B|<\tet$, then we can shrink $S$  to set
$S'$  of the same cardinality such that for
every $\bet, \alp\in S'$  $\onu_{\alp
+1}=\onu_{\bet +1}$.  Now projecting down to limit
points of $S'$ of uncountable cofinality we will
obtain (*) of the conclusion of the theorem.
So, suppose now that $|B|=\tet$.  W.l. of g., we
can assume that $\alp <\bet$  implies
$\onu_{\alp +1}<\onu_{\bet +1}$.  Now, by [Git-Mit],
$B$ (or at least its initial segments) is
contained in the length of an extender over $\nu_\tet$
in the core model.  There is no overlapping
extenders, hence 
$$tcf\Big(\prodl_{\alp\in S}\onu_{\alp +1}\Big/J^{bd}_\tet
\Big)=\Big(\sup (\{\onu_{\alp+1}\mid\alp\in S\})\Big)^+$$
where the successor is in sense of the core model
or the universe which is the same by the Mitchell
Weak Covering Lemma.  Also, for every $\alp$
which is a limit point of $S$  of uncountable cofinality
$$tcf\Big(\prodl_{\bet\in S\cap_\alp}\onu_{\bet +1}\Big/
J^{bd}_{S\cap\alp}\Big)=\Big(\sup \{\onu_{\bet
+1}\mid\bet\in S\cap\alp\}\Big)^+\ .$$
Projecting down we obtain (*). 

\medskip
\noindent
{\bf Subcase 2.2.}\quad Starting with some $\alp^* <\tet$
each $\nu^*_{\alp +1}$ is not a principal
indiscernible for $\nu_\tet$  or it is but
$\nu_{\alp +1}$  corresponds over $\nu_\tet$
to some $\onu_{\alp +1}$  which has cofinality
$<\nu_\tet$  in the core model.

Suppose for simplicity that $\alp^*=0$.  If
$\nu^*_{\alp +1}$  is not a principal indiscernible for 
$\nu_{\tet}$, then we can use functions of the
core model to transfer the structure of
indiscernibles over $\nu^*_{\alp +1}$ to the
interval $[\nu_\alp$, length of the extender used
over $\nu_\alp]$. This will replace $\nu_{\alp
+1}$ be a member of the interval.  So let us
concentrate on the situation when $\nu^*_{\alp
+1}$  is a principal indiscernible for $\nu_\tet$  
but $\onu_{\alp +1}$  has cofinality $\le\nu_\tet$
$(\alp <\tet)$.

Let us argue that this situation is impossible.
Thus we have increasing sequences
$\langle\alp_i\mid i\le\tet\rangle$, $\langle
\rho_i\mid i <\tet\rangle$ and $\langle\rho'_i\mid
i<\tet\rangle$  such that for every $i<\tet$
$\rho_i$  is between $\nu_{\alp_i}$  and the
length of the extender used over $\nu_{\alp_i}$,
$cf\rho_i\ge\nu_{\alp_i}$  in the core model,
$\rho'_i$  is the image of $\rho_i$  over
$\nu_{\alp_i+1}$  and $cf\rho'_i<\nu_{\alp_i+1}$  in the
core model.  Then $cf\rho'_i<\nu_{\alp_i}$  again in the
core model since
$\rho'_i$ is the image of $\rho_i$  in the
ultrapower and $\nu_{\alp_i+1}$  the image of
$\nu_{\alp_i}$  which is the critical point of
the embedding.  Fix for every $i<\tet$  a
sequence $c_i$  unbounded in $\rho'_i$,  in the
core model and of cardinality $cf\rho'_i$ there.
Take a precovering set including $\{c_i\mid
i<\tet\}$.  By [Git-Mit], assignment functions
can change for this new precovering set only on
a bounded subset of $\nu_{\alp_i}$'s.  Pick
$i<\tet$  such that $\nu_{\alp_i}$  is above
supremum of this set.  Again, consider the
ultrapower used to move from $\nu_{\alp_i}$  to
$\nu_{\alp_i+1}$. Now we have $c_i$  in this
ultrapower and its cardinality is $<\nu_{\alp_i}$.
Let $j:M\to M'$ be the embedding.  $c_i\in M'$
and $M'$  is an ultrapower by extender.  Hence
for some $\tau$  and $f$ $c_i=j(f)(\tau)$.  Let
$U_\tau=\{X\subseteq\nu_{\alp_i}\mid\tau\in j(X)\}$  and  
$\tilj:M\to\tilM$  be the corresponding
ultrapower.  Denote $\tilj(\nu_{\alp_i})$ by
$\tilnu_{a_i +1}$,  $\tilj(\rho_i)=\tilc_i$
and $\tilj(f)([id])=\tilc_i$.  Let $\tilc_i=\langle
j(f_\xi)([id])\mid \xi <\xi^*=cf\rho'_i=cf\tilrho_i
\rangle$ be increasing enumeration (everything
in the core model).  Then for most $\bet$'s
$(\MOD U_\tau)$ $f(\bet)=\langle f_\xi(\bet)\mid
\xi <\xi^*\rangle$ will be a sequence in $M$
cofinal in $\rho_i$  of order type $\xi$.
Which contradicts the assumption that
$cf\rho_i\ge\nu_i$.\hfill $\square$

Let us use 1.9 in order to deduce the following:

\medskip
\noindent
{\bf Theorem 1.10.}\quad {\sl Suppose that there is no
inner model with strong cardinal then for every set $a$
of regular cardinals above
$2^{|a|^++\aleph_2}$ $\ \ |pcfa|\le |a|+\aleph_1$. 
}

\medskip
\noindent
{\bf Remark.}\quad If $a$  is an interval then
$|pcf a|=|a|$ by [Git-Mit, 3.24].

\medskip
\noindent
{\bf Proof.}\quad Suppose that for some $a$  as
in the statement of the theorem $|pcf
a|>|a|+\aleph_1$.
Let $\tet =|a|^++\aleph_2$.  Then $|pcf a|\ge
\tet$.  Pick an increasing sequence
$\langle\tau_\alp\mid\alp <\tet\rangle$  inside
$pcf (a)$.  By 1.9 we can find an unbounded
subset $S$  of $\tet$  satisfying the conclusion
(*) of 1.9.

Let $D$  be an ultrafilter on $\tet$  including
all cobounded subsets of $S$.  Let $\tau =
cf(\prodl_{\alp <\tet}\tau_\alp/D)$.  Then, clearly,
$\tau\ge(\bigcup_{\alp <\tet}\tau_\alp)^+$.  By the
Localization Theorem [Sh-g], then there is
$a_0\subseteq\{ \tau_\alp\mid\alp\in S\}$, $|a_0|\le |a|$
with $\tau\in pcf a_0$.  Consider $S\bks\sup a_0$.
$S\bks\sup a_0\in D$ since $a_0$ is bounded in $S$.  Hence
$cf\Big(\prodl_{\alp\in S\bks\sup a_0}\tau_\alp\Big/D
\Big)=\tau$.  Again by the Localization Theorem,
there is $a_1\subseteq S\bks\sup a_0$, $|a_1|\le
|a|$  and $\tau\in pcf a_1$.  Continue by
induction and define a sequence $\langle a_\alp
\mid \alp <\ome_1\rangle$  such that for every
$\alp <\ome_1$ the following holds:
\begin{itemize}
\item [{\rm (a)}] $a_\alp\subseteq S$
\smallskip
\item [{\rm (b)}] $|a_\alp|\le |a|$
\smallskip
\item [{\rm (c)}] $\tau\in pcf a_\alp$
\smallskip
\item [{\rm (d)}] $\min a_\alp >\sup a_\bet$  for every
$\bet<\alp$.
\end{itemize}

Let $\del =\bigcup_{\alp <\ome_1}\sup a_\alp$.
Then $\del$  is a limit of points of $S$  and $cf
\del =\aleph_1$.  Hence (*) of 1.9 applies.
Thus $tcf\Big(\prodl_{\alp\in\del\cap
S}\tau_\alp / J^{bd}_{\del\cap S}\Big)$ exists
is below $\tau_{\del +1}$  and is equal to
$tcf\Big(\prodl_{\alp\in\del\cap S}\tau_\alp/F\Big)$  
for every ultrafilter $F$  on $\del\cap S$
including all cobounded subsets of $\del\cap
S$.  Denote  $tcf\Big(\prodl_{\alp\in\del\cap S}
\tau_\alp/J^{bd}_{\del\cap S}\Big)$  by $\mu$.
Let $c=pcf (a)$  and $\langle b_\xi[c]\mid\xi\in
pcf (a)=c\r$ be a generating sequence. Clearly
both $\mu$  and $\tau$  are in $c$  and $\mu <\tau$.    
Consider $b=b_\tau [c]\bks b_\mu [c]$.

For every $\alp <\ome_1$, $b\cap a_\alp\not=\emptyset$,
since $\tau\in pcf (a_\alp)$. Hence, $b\cap\del\cap S$
is unbounded in $\del$ (by (d) of the choice
of $a_\alp$'s).  Let $F$  be an ultrafilter on
$\del\cap S$  including $b\cap\del\cap S$ And
all cobounded subsets of $\del \cap S$.  Then
$tcf \Big(\prodl_{\alp\in\del\cap S}\tau_\alp/
F\Big) =\mu$ but  this means that $\mu\in pcf b$, which
is impossible by the choice of $b$,  see for
example [Bur-Mag, 1.2].\hfill $\square$    

The proof of 1.10 easily gives a result related
to the strength of the negation of the Shelah
Weak Hypothesis (SWH). 
(SWH says that for every cardinal $\lam$ the
number of singular cardinals $\kap <\lam$  with
$pp\kap\ge\lam$  is at most countable).

\medskip
\noindent
{\bf Theorem 1.10.1.}\quad {\sl Suppose that
there is no inner model with strong cardinal.
Then for every cardinal $\lam >2^{\aleph_2}$ 
$$|\{\kap <\lam \mid cf\kap <\kap\quad\text{and}\quad
pp\kap\ge\lam \}|\le \aleph_1\ .$$
}

Now we continue the task started in 1.1. and
deal with ordinal gaps.

Let us start with technical definitions.

\medskip
\noindent
{\bf Definition 1.11.}\quad Let 
$$\text{Kinds}=\Big\{\del^{\ell_0}_0\cdot\del^{\ell_1}_1
\cdots\del^{\ell_{k-1}}_{k-1}\Big| k <\ome,
1\le\ell_0\nek \ell_{k-1} <\ome, \del_0 >\del_1
>\cdots\del_{k-1}\quad\text{are cardinals}$$
of uncountable cofinality $\Big\}\cup\{ 0\}$, where the
operations used are the ordinals operations.

\medskip
\noindent
{\bf Remark 1.12.}\quad The only kinds around
$\ome_1$  are $\ome_1$  itself, $\ome^2_1\nek\ome^n_1
\cdots (n <\ome)$.  But already with $\ome_2$
we can generate in addition to $\ome_2,\ome^2_2\nek
\ome^n_2\cdots (n<\ome)$  also $\ome_2\cdot\ome^5_1$,
$\ome^{19}_2\cdot\ome^3_1$ etc.  Note that between
$\ome^\ome_1$  and $\ome_2$ there are no new kinds.
Using the Rado-Milner paradox we will show in the next
section that the consistency strength of the length the
gap does not change in such an interval.  

\medskip
\noindent
{\bf Definition 1.13.}\quad Let $\gam$  be an ordinal
\begin{itemize}
\item[{\rm (a)}] $\gam$  is of kind $0$  if
$\gam$ is a limit ordinal.
\item[{\rm (b)}] $\gam$  is of kind  $\del_0$
for a cardinal $\del_0\in {\rm Kinds}$ if $\gam$ is
a limit of an increasing sequence of length
$\del_0$.  In particular, if $\del_0$  is regular
this means that $cf\gam =\del_0$. 
\item[{\rm (c)}] $\gam$  is of kind
$\del_0^{\ell_0}\cdot \del_1^{\ell_1}\cdots
\del^{\ell_{k-1}}_{k-1}\in {\rm Kinds}$, with
$\ell_0>1$, if $\gam$ is a limit of an
increasing sequence of $\del_{k-1}$  ordinals of
kind $\del^{\ell_0}_0\cdot\del^{\ell_1}_1\cdots  
\del^{\ell_{k-1}-1}_{k-1}$.
\end{itemize}

\medskip
\noindent
{\bf Lemma 1.14.}\quad {\sl Let $\kap$  be a strong limit
cardinal of cofinality $\aleph_0$, $\del <\kap$ a
cardinal of uncountable cofinality.  Assume
\begin{itemize}
\item [(1)] $SSH^{\le\del}_{<\kap}$
\item [(2)] there is no measurable cardinals in
the core model between $\kap$ and
$\kap^{+\del^+}$.
\end{itemize}

Let $0<\xi\in {\rm Kinds}\cap [\del,\del^+)$
and $2^\kap\ge\kap^{+\alp +\xi}$, for some $\alp
<\del^+$.  Then $\kap^{+\alp +\xi +1}\in pcf
\{\tau^{+\nu +1}_{n,i}\mid\nu\quad\text{is an
ordinal of kind}\quad \xi, i <i(n), n<\ome\}$,
where $\tau_{ni}$ denotes the principal
indiscernible of the block $B_{n,i}$,  as
defined in 1.6. 
}

\medskip
\noindent
{\bf Remark 1.15.}
\begin{itemize}
\item [{\rm (a)}] The lemma provides a bit more
information then will be needed for deducing the
strength of $2^\kap =\kap^{+\xi +1}$.  
\item [{\rm (b)}] The condition (2) is not very
restrictive since we are interested in small
$(<\kap)$  gaps between $\kap$  and its power. 
\end{itemize}

\medskip
\noindent
{\bf Proof.}\quad We prove the statement by induction
on $\xi$.  Fix $\alp <\del^+$.  Let $\xi
=\del^{\ell_0}_0\cdots\del^{\ell_{k-1}}_{k-1}$,
where $\del_0=\del$.  Set for each $\sig <\del_{k-1}$
$$\kap(\sig)=\kap^{+\alp +\del_0^{\ell_0}\cdots
\del^{\ell_{n-2}}\cdot\del^{\ell_{n-1}-1}_{k-1}\cdot\sig
+\del_0^{\ell_0}\cdots\del^{\ell_{k-2}}_{k-2}\cdot
\del^{\ell_{k-1}-1}_{k-1}+1}$$
if $(k>1)$ or $(k=1$  and $\ell_0 >1$) and
$$\kap (\sig)=\kap^{+\alp +\sig +1}$$
if $k=1$  and $\ell_0=1$, i.e. $\xi=\del$.

For every $\sig <\del_{k-1}$, if $\xi\not=\del$
then by induction
\newline
$\kap(\sig)\in pcf(\{\tau^{+\nu +1}_{n,i}|\nu$  is an
ordinal of kind $\del_0^{\ell_0}\cdots
\del_{k-2}^{\ell_{k-2}}\cdot\del_{k-1}^{\ell_{k-1}-1},
i<i(n), n<\ome\})$.

Let $E$  be the set consisting of all regular
cardinals of blocks $B_{n,i} (n<\ome, i<i(n))$
together with all regular cardinals between
$\kap$  and $\min\Big(\kap^{+\del^+} ,2^\kap)$.
Set $E^*=pcf E$. Then $\kap >|pcf E^*|$,  since
$\kap$  is strong limit.  We can assume also
that $\min E^*>|pcf E^*|$.  By [Sh-g], then $pcf
E^*=E^*$ and there is a set $\langle b_\chi
[E^*]\mid\chi\in E^*\rangle$  of $pcf E^*$
generators which is smooth and closed, i.e.
$\tau\in b_\chi [E^*]$ implies $b_\tau [E^*]\subseteq  
b_\chi [E^*]$  and $pcf (b_\chi [E^*])=b_\chi
[E^*]$. 

The assumption (2) of the lemma implies that for
every unbounded in $\kap^{+\alp+\xi}$ set $B$
consisting of regular cardinals above $\kap$ and
below $\kap^{+\alp +\xi}$ $\max pcf(B)=\kap^{+\alp
+\xi +1}$.  In particular $\max pcf \{\kap
(\sig)|\sig <\del_{k-1}\})=\kap^{+\alp +\xi +1}$.  Denote
$\kap^{+\alp +\xi +1}$  by $\mu$.  Let 
$$A^*=b_\mu [E^*]\cap \{\kap (\sig)\mid\sig
<\del_{k-1}\}\ .$$
Then, $|A^*|=\del_{k-1}$ and for every $\lam\in
A^*$ $b_\lam [E^*]\subseteq b_\mu[E^*]$.  For
every $\lam\in A^*$, fix a sequence $\langle
\rho^\lam_n\mid
n<\ome\rangle\in\prodl_{n<\ome}\kap^+_{n+1}$
inside $b_\lam[E^*]$ such that
\begin{itemize}
\item [({\rm a})] $\rho^\lam_n\in B_{n,i}$  for
some $i<i(n)$
\newline
and, if $\xi\not=\del$ then also 
\item [({\rm b})] $\rho^\lam_n$  is of kind
$\del_0^{\ell_0}\cdots\del_{k-2}^{\ell_{k-2}}\cdot
\del_{k-1}^{\ell_{k-1}-1}$.  
\end{itemize}
It is possible to find $\rho^\lam_n$'s of the
right kind using the inductive assumption, as
was observed above.

\medskip
\noindent
{\bf Claim 1.16.}\quad There are infinitely many
$n<\ome$ such that 
$$|\{\rho^\lam_n\mid\lam\in a^*\}|=\del_{k-1}$$

\medskip
\noindent
{\bf Proof.}\quad Otherwise by removing finitely
many $n$'s or boundedly many $\rho^\lam_n$'s we
can assume that for every $n$ $|\{\rho^\lam_n\mid\lam\in 
A^*\}|<\del_{k-1}$.  But $cf\del_{k-1} >\aleph_0$.  Hence,
the total number of $\rho^\lam_n$'s is less than
$\del_{k-1}$.  Now, $pcf\{\rho^\lam_n\mid n<\ome\ ,\ 
\lam\in A^*\}\supseteq A^*$.  So, $|A^*\cap
pcf\{\rho^\lam_n\mid n<\ome,\lam\in A^*\}|\ge
|A^*|=\del_{k-1}$.  By (2) of the statement of
the lemma this situation is impossible.

\hfill $\square$ of the claim.

Suppose for simplicity that each $n<\ome$
satisfies the conclusion of the claim.  If not
then we just can remove all the ``bad" $n$'s.
This will effect less than $\del_{k-1}$  of $\rho$'s
which in turn effects less than $\del_{k-1}$ of $\lam$'s.

Let us call a cardinal $\tau$ reasonable, if for
some $n<\ome$ $\tau$ is a limit of $\del_{k-1}$-sequence
of elements of $\{\rho^\lam_n\mid \lam\in A^*\}$.
Clearly, a reasonable $\tau$  is of kind
$\del_0^{\ell_0}\cdots \del_{k-1}^{\ell_{k-1}}$,
since $\rho^\lam_n$'s are of kind
$\del_0^{\ell_0}\cdot\del_1^{\ell_1}\cdots
\del_{k-2}^{\ell_{k-2}}\cdot\del_{k-1}^{\ell_{k-1}-1}$.
The successor of such $\tau$ is in $pcf\{\rho^\lam_n
\mid\lam\in A^*\}$  since $cf\tau=cf\del_{k-1}$  and we
assumed $SSH^{cf\del_{k-1}}_{<\kap}$,
i.e.  $pp\tau =\tau^+$.  Also $pp\tau=\tau^+$
implies that the set $\{\rho^\lam_n\mid\lam\in
A^*\}\bks b_{\tau^+}[E^*]$  is bounded in $\tau$.

\medskip
\noindent
{\bf Claim 1.17.}\quad $pcf\{\tau^+\mid\tau$  is
reasonable$\}\subseteq b_\mu[E^*]$. 

\medskip
\noindent
{\bf Proof.}\quad $\{\rho^\lam_n\mid n<\ome\}\subseteq
b_\lam [E^*]$ for every $\lam\in A^*$.  Also,
$b_\lam [E^*]\subseteq b_\mu [E^*]$.  By the above,
for every reasonable $\tau$, $\tau^+\in
pcf\{\rho^\lam_n\mid\lam\in A^*\}$  for some
$n<\ome$.  But $pcf(b_\mu[E^*])=b_\mu [E^*]$ and
$pcf\{\rho^\lam_n| n<\ome,\lam\in A^*\}\subseteq
pcf (b_\mu [E^*])$ since the pcf generators are
closed and $\{\rho^\lam_n\mid n<\ome,\lam\in
A^*\}\subseteq b_\mu[E^*]$.  So,
$\{\tau^+\mid\tau\ \hbox{is reasonable}\}\subseteq
b_\mu [E^*]$ and again using closedness of $b_\mu[E^*]$,
we obtain the desired conclusion.  

\hfill $\square$  of the claim.

\medskip
\noindent
{\bf Claim 1.18.}\quad For every
$\mu'\in pcf\{\tau^+\mid\tau\ \hbox{is reasonable}\}$, 
$b_{\mu'}[E^*]\subseteq b_\mu[E^*]$.

\medskip
\noindent
{\bf Proof.}\quad By the smoothness of the generators
$b_{\mu'}[E^*]\subseteq b_\mu[E]$  for every
$\mu'\in pcf\{\tau^+\mid\tau\ \hbox{is
reasonable}\}$.
\newline
\hfill $\square$  of the claim.

In order to conclude the proof we shall argue
that there should be $\mu'\in pcf\{\tau^+\mid\tau$
is reasonable$\}$  such that $\mu\in
b_{\mu'}[E^*]$.  This will imply $b_\mu
[E^*]=b_{\mu'}[E^*]$ and hence $\mu=\mu'$.

Let us start with the following:

\medskip
\noindent
{\bf Claim 1.19.}\quad 
$|\{\rho^\lam_n\mid n<\ome,\lam\in A^*\}\bks\bigcup
\{b_{\tau^+}[E^*]|\tau\ \hbox{is
reasonable}\}|<\del_{k-1}$. 

\medskip
\noindent
{\bf Proof.}\quad Suppose otherwise.  Let
$S=\{\rho^\lam_n\mid n<\ome,\lam\in A^*\}\bks
\bigcup\{b_{\tau^+}[E^*]|\tau\ \hbox{is reasonable}\}$
and $|S|=\del_{k-1}$.  Then for some $n<\ome$  also
$\{\rho^\lam_n\mid\rho^\lam_n\in S\}$  has
cardinality $\del_{k-1}$, since $cf\del_{k-1} >\aleph_0$.
Fix such an $n$  and denote
$\{\rho^\lam_n\mid\rho^\lam_n\in S\}$  by $S_n$.

But now there is a reasonable $\tau$  which is a
limit of elements of $S_n$.  $pp\tau =\tau^+$
implies that the set $\{\rho^\lam_n\mid\lam\in
A^*\}\bks b_{\tau^+}[E^*]$ is bounded in $\tau$.
In particular, $S_n\cap b_{\tau^+}[E^*]$  is
unbounded.  Contradiction, since $S_n\subseteq
S$  which is disjoint to every $b_{\tau^+}[E^*]$
with $\tau$ reasonable.

\hfill $\square$ of the claim.

Now, removing if necessary less than $\del$
elements, we can assume that $\{\rho^\lam_n\mid
n<\ome,\lam\in A^*\}$  is contained in $\cup\{b_{\tau^+}
[E^*]\mid\tau\ \hbox{is reasonable}\}$.  Recall
that this can effect only less than $\del$  of
$\lam$'s in $A^*$  which has no influence on
$\mu$.

Let $b=pcf\{\tau^+\mid\tau\ \hbox{is reasonable}\}$. 
Then $pcf b=b$  and $b\subseteq E^*$.  By [Sh-g], there
are $\mu_1\nek\mu_\ell\in pcf b=b$ such that
$b\subseteq b_{\mu_1}[E^*]\cup\cdots\cup b_{\mu_\ell} 
[E^*]$.  Using the smoothness of generators, we
obtain that for every reasonable $\tau$  there
is $k$, $1\le k\le\ell$  such that
$b_{\tau^+}[E^*]\subseteq b_{\mu_k}[E^*]$. Now,
$\{\rho^\lam_n\mid n<\ome,\lam\in A^*\}\subseteq\cup
\{b_{\tau^+}[E^*]\mid\tau\ \hbox{is reasonable}\}$.
Hence, $\{\rho^\lam_n\mid n<\ome,\lam\in
A^*\}\subseteq\bigcup^\ell_{k=1}b_{\mu_k}[E^*]$.

For every $\lam\in A^*$  fix an ultrafilter
$D_\lam$  on $\ome$  including all cofinite sets
so that $tcf\Big(\prodl_{n<\ome}\rho^\lam_n\Big/D_\lam
\Big)=\lam$.  Let $\lam\in A^*$.  There are
$x_\lam\in D_\lam$  and $k(\lam)$, $1\le
k(\lam)\le\ell$  such that for every $n\in
x_\lam$ $\rho^\lam_n\in b_{\mu_{k(\lam)}}[E^*]$.
Then $\lam\in pcf\Big(b_{\mu_{k(\lam)}}[E^*]\Big)=
b_{\mu_{k(\lam)}}[E^*]$.  Finally, we find
$A^{**}\subseteq A^*$ of cardinality $\del_{k-1}$ (or
just unbounded in $\mu$) and $k^*$, $1\le
k^*\le\ell$  such that for every $\lam\in
A^{**}$  $k(\lam)=k^*$.  Then $A^{**}\subseteq
b_{\mu_{k^*}}[E^*]$.  But, recall that $\mu =\max pcf (B)$
for every unbounded subset $B$ of $A^*$.  In
particular, $\mu =\max pcf (A^{**})$.  Hence,
$\mu\in pcf A^{**}\subseteq pcf\Big(b_{\mu_{k^*}}[E^*]
\Big)=b_{\mu_k^*}[E^*]$.
\newline
\hfill $\square$

Lemma 1.14 implies the following:

\medskip
\noindent
{\bf Theorem 1.20.}\quad {\sl Let $\kap$  be a strong limit
cardinal of confinality $\aleph_0$, $0<\xi\in {\rm Kinds}$.
Assume that 
\begin{itemize}
\item [(1)] $SSH^{\le |\xi |}_{<\kap}$
\item [(2)] there are no measurable cardinals in the core model
between $\kap$ and $\kap^{+|\xi|^+}$.
\end{itemize}

If $2^\kap\ge\kap^{+\xi}$, then in the core model either
\begin{itemize}
\item [(i)] $o(\kap)\ge\kap^{+\xi +1}+1$
or
\item [(ii)] $\{\alp <\kap\mid o(\alp)\ge\alp^{+\xi +1}+1\}$ is
unbounded in $\kap$.
\end{itemize}
}

\medskip
\noindent
{\bf Proof.}\quad  By 1.14, for infinitely many $n$'s for some
$i_k <i(n)$  the length of the block $B_{n,i*_n}$  will be at
least $\tau^{+\xi +1}_{n,i_n}$, since it should contain some
$\tau^{+\nu +1}_{n,i_n}$ for an ordinal $\nu$  of kind $\xi$.
Clearly, $\nu\ge\xi$  since $\xi$  is the least ordinal of kind
$\xi$.
\newline
\hfill $\square$

We like now outline a way to remove (2) of 1.20 by cost of
restricting possible $\xi$'s.
First change Definitions 1.11 and 1.13.  Thus in 1.11 we
replace uncountable by ``above $\aleph_1$".  Denote
by Kinds$^*$ the resulting class.  Then define kind$^*$ of
ordinal as in 1.13 replacing Kinds by Kinds$^*$.

\medskip
\noindent
{\bf Theorem 1.21.}\quad {\sl Let $\kap$  be a strong
limit cardinal of cofinality $\aleph_0$,
$0<\xi\in {\rm Kinds}^*$.  Assume
$SSH^{\le |\xi |}_{<\kap}$.  If $2^\kap\ge\kap^{+\xi}$,
 then in the core model either
\begin{itemize}
\item [(i)] $o(\kap)\ge\kap^{+\xi +1}+1$
or
\item [(ii)] $\{\alp<\kap\mid o(\alp)\ge\alp^{+\xi+1}+1\}$
is unbounded in $\kap$.
\end{itemize}
}
The theorem, as in the case of 1.20, will follow from the
following: 

\medskip
\noindent
{\bf Lemma 1.22.}\quad {\sl Let $\kap$ be a strong limit cardinal
of cofinality $\aleph_0$, $\del <\kap$ a cardinal of cofinality
above $\aleph_1$.  Assume $SSH^{\le\del}_{<\kap}$.  Let $0 <\xi\in
{\rm Kinds}^*\cap [\del,\del^+]$ and $2^\kap\ge\kap^{+\alp +\xi}$
for some $\alp <\del^+$.  Then
\newline
$$pcf\{\tau^{+\nu +1}_{ni}\mid
\nu\quad\text{is an ordinal of kind}^*\ \xi,i<i(n),n
<\ome\}\cap [\kap^{+\alp+\xi+1},\kap^{+\alp +\xi+\xi +1}]
\not=\emptyset\ ,$$
}

Let us first deal with a special case -- $\xi$  is
a cardinal.  We
split it into two cases: (a) $\xi$  is regular and (b) $\xi$ is
singular.  The result will be stronger than those of 1.22.  

\medskip
\noindent
{\bf Lemma 1.23.}\quad {\sl Let $\kap$ be a strong limit cardinal
of cofinality $\aleph_0$, $\del <\kap$ is a regular uncountable
cardinal.  Assume $SSH^{\le\del}_{<\kap}$.  Let
$2^\kap\ge\kap^{+\alp +\del}$ for some $\alp >\del^+$.  Then
$$\kap^{+\alp +\del +1}\in pcf\Big(\Big\{\tau^{+\nu +1}_{ni}\mid
i<i(n),\ n <\ome\quad\text{and}\quad\nu\quad\text{is an ordinal
of cofinality}\quad \del\}\Big)\ .$$
}

\medskip
\noindent
{\bf Proof.}\quad Let $\mu =\kap^{+\alp +\del +1}$.
 We choose $E^*$ and $\langle b_\chi[E^*]\mid\chi\in
E^*\rangle$ as in the proof of 1.14.  Measurables of a
core model between $\kap$  and $2^\kap$  are allowed
here.  So in contrast to 1.14 we cannot
claim anymore for every unbounded
$B\subseteq [\kap,\kap^{+\alp +\del})$ consisting of
regulars $\max pcf (B)=\kap^{+\alp +\del +1}$.
Hence the choice of $A^*$  (the crucial for the proof
set in 1.14) will be more careful.       

Set $A$  to be the set of cardinals
$\kap^{+\alp +\tau +1}\in [\kap^{+\alp +1},
\kap^{+\alp+\del})$  such that either $o(\bet)<
\kap^{+\alp +\tau}$ for every $\bet <\kap^{+\alp +\tau}$
or else $\kap^{+\alp +\tau}$  is above every measurable
of the core model smaller than $\kap^{+\alp +\del}$.
Clearly, $|A|=\del$, since there is no overlapping
extenders and as in 1.1 $|(pcf
b)\cap A|\le |b|$  for every set of regular cardinals
$b\subseteq\kap$, $|b|\le\del$.  By 1.3, $\max pcf (B)=
\kap^{+\alp +\del +1}$  for every unbounded
$B\subseteq A$.  This implies that $A\bks b_\mu [E^*]$
 is bounded in $\kap^{+\alp +\del +1}$.  Define
$A^*=A\cap b_\mu [E^*]$. The rest of the proof
completely repeats 1.14.
\newline
\hfill $\square$

\medskip
\noindent
{\bf Lemma 1.24.}\quad {\sl Let $\kap$  be a strong limit
cardinal of cofinality $\aleph_0$, $\del <\kap$  is a
regular cardinal of uncountable cofinality.  Assume
$SSH^{\le\del}_{<\kap}$.
Let $2^\kap\ge\kap^{+\alp +\del}$  for some
$\alp <\del^+$.  Then
\newline
$pcf (\tau^{\nu +1}_{ni}|i<i(n),n<\ome$ and $\nu$  is a
limit of an increasing sequence of the length $\del\})
\cap [\kap^{+\alp +\del +1},\kap^{+\alp +\del +\del +1}]
\not=\emptyset$.
}

\medskip
\noindent
{\bf Proof.}\quad Let $\langle\del_i\mid i<cf\del\rangle$  be an
increasing continuous sequence of limit cardinals unbounded in
$\del$.  Consider the set
$$B=\{\kap^{+\alp +\del_i+\nu}\mid i<cf\del\ ,\
i\quad\hbox{limit and}\quad \nu <\del_i\} .$$ 
Since $cf\del >\aleph_0$, the analysis of indiscernibles of
[Git-Mit, Sec. 3.4] can be applied to show that $\{cf (\prod
B/D)|D$  is an ultrafilter over $B$ extending the filter of
cobounded subsets of $B\}\subseteq\{\kap^{+\alp +\nu +1}\mid\del  
\le\nu\le\del +\del \}$. 

We cannot just stick to $\kap^{+\alp +\del +1}$ alone since we
like to have $\del$  cardinals below $\kap^{+\alp +\del}$.  But
once measurable above $\kap$ allowed, it is possible that $\max
pcf(\{\kap^{+\alp +\rho +1}\mid\rho <\del \})>\kap^{+\alp +\del +1}$.
Still by [Sh-g], for a club $C\subseteq cf\del\quad tcf
(\prodl_{\nu\in C}\kap^{+\alp +\nu +1}/$ cobounded $\upharpoonright
C)=\kap^{+\alp+\del +1}$.  Unfortunately, this provided only
$cf\del$ many cardinals $\kap^{+\alp +\nu +1}$  and not $\del$-many.

Define a filter $D$  over $B$:
\newline
 $X\in D$  iff $\{i <cf\del|i$
is limit and $\{ j<i|\{\nu <\del^+_j|\kap^{\alp +\del_i+\xi
+1}\in X\}$  is cobounded in $\del^+_j\}$  is cobounded in $i\}$
contains a club.

Let $D^*$  be an ultrafilter extending $D$.  Set $\mu =cf (\prod
B/D^*)$.  By the choice of $D$, for every $C\subseteq B$  of cardinality
less than $\del$ $B\bks C\in D$.  So, $\mu\in [\kap^{+\alp +\del
+1}, \kap^{+\alp +\del +\del +1}]$.  Define $E^*$  as before.
Set $A^*=B\cap b_\mu[E^*]$.

\medskip
\noindent
{\bf Claim 1.25.}\quad If $A^*\in D^*$.

\medskip
\noindent
{\bf Proof.}\quad Otherwise the compliment of $A^*$  is in $D^*$.
Let $A'=B\bks b_\mu[E^*]$.  Clearly, $D^*\cap J_{<\mu}[E^*]
=\emptyset$.  By [Bur-Mag, 1.2], then there is $S\in D^*$ $S\in
J_{<\mu^*}[E^*]\bks J_{<\mu}[E^*]$. But $b_\mu [E^*]$  generates
$J_{<\mu}[E^*]$ over $J_{<\mu}[E^*]$.  So, $S\subseteq b_\mu[E^*]
\cup c$  for some $c\in J_{<\mu}[E^*]$.  Hence, $S\cap b_\mu
[E^*]\in D^*$.  But $A'\in D^*$ and $A'\cap B\cap (S\cap b_\mu
[E^*])=\emptyset$.  Contradiction.
\newline
\hfill $\square$ of the claim.

Now we continue as in the proof of 1.14.  In order to eliminate 
possible effects of less than $\del$  cardinals, we use 1.10.  At
the final stage of the proof a set $A^{**}$  was defined.  Here
we pick it to be in $D^*$.  This insures that $\mu\in pcf A^{**}$
and we are done.
\newline
\hfill $\square$

Now we turn to the proof of 1.22.

\medskip
\noindent
{\bf Proof.}\quad As in 1.14, we prove the statement by induction
on $\xi$. Fix $\alp <\del^+$.  Let $\xi =\del_0^{\ell_0}\cdots
\del_{k-1}^{\ell_{k-1}}$. The case $k=1$  and $\ell_0=1$  (i.e.
$\xi =\del$) was proved in 1.23, 1.24. So assume that $k>1$  or
$(k=1$  and $\ell_0 >1)$.  For each $\sig <\del_{k-1}$ let 
$$\kap (\sig)\in pcf (\{\tau^{+\nu +1}_{n,i}\mid i<i(n),
n<\ome\ ,\quad\text{and}\quad\nu\quad\text{is an ordinal
of kind}^*\quad \del_0^{\ell_0}\cdots
\del_{k-2}^{\ell_{k-2}}\cdot\del_{k-1}^{\ell_{k-1}-1}\})
\cap$$
$$[\kap^{+\alp+\xi^-\cdot\sig +\xi^-+1},\kap^{+\alp
+\xi^-\cdot\sig+ \xi^-+ \xi^-+1}]\ ,\quad {\rm where}$$
$$\xi^-=
\begin{cases}
\del_0^{\ell_0}\cdots\del_{k-2}^{\ell_{k-2}}\cdot
\del_{k-1}^{\ell_{k-1}-1}\ ,if&\xi =\del_0^{\ell_0}
\cdots\del_{k-1}^{\ell_{k-1}}\quad\text{and}\quad
(k>1\quad {\rm or}\quad (k=1\quad {\rm and}\quad
\ell_0>1)\\
0,&if\quad k=1\quad\text{and}\quad\ell_0=1
\end{cases}
$$
In the last case the inductive assumption insures the existence
of such $\kap (\sig)$.

Define $E^*$  and $\langle b_\chi[E^*]|\chi\in E^*\rangle$ as in
the proof of 1.14. We do not know now if for every unbounded in
$\kap^{+\alp +\xi}$  set $B\subseteq [\kap,\kap^{+\alp +\xi})$
consisting of regular cardinals $\max pcf (B)=\kap^{+\alp +\xi
+1}$.  We may consider the set $\{\kap^{+\alp +\xi^-\cdot\nu
+1}\mid\nu <\del_{k-1}\}$.  If for club many $\nu$'s
$\kap^{+\alp +\xi^-\cdot\nu +1}$ is not a principle indiscernible
then by [Git-Mit] $cf (\prod B/{\rm bounded})=\kap^{+\alp +\xi +1}$
for any unbounded subset $B$  of $\kap^{+\alp +\xi}$  consisting
of regular cardinals.  Note that $cf\del_{k-1}>\aleph_0$  is
crucial here.  In this case we define $A^*=\{\kap(\sig)\mid\sig
<\del_{k-1}\}\cap b_{\kap^{+\alp +\xi +1}}[E^*]$  and  proceed as
in the proof of 1.14.  The only difference will be the use of
1.10 to eliminate a possible influence of $<\del_{k-1}$ cardinals.
Here the assumption $\del_{k-1}>\aleph_1$  comes into play.
In the general case it is possible to have $\{\kap (\sig )\mid
\sig <\del_{k-1}\}\cap b_{k^{+\alp+\xi+1}}[E^*]$  empty.  But
once for a club of $\nu$'s below $\del_{k-1}$ $\kap^{\alp
+\xi^-\cdot\nu +1}$'s are principal indiscernibles, by [Git-Mit]
we can deduce that
$$pcf(\{\kap (\sig)\mid\sig <\del_{k-1}\})\bks
\kap^{+\alp+\xi}\subseteq$$
$$[\kap^{+\alp+\xi +1},\kap^{+\alp +\xi +\xi^-+\xi^-+1}]\subseteq
[\kap^{+\alp +\xi +1},\kap^{+\alp +\xi +\xi +1}]\ .$$

Let $D$  be an ultrafilter on the set $\{\kap (\sig)\mid\sig
<\del_{k+1}\}$  containing all cobounded subsets.  Set
$$\mu =cf (\prod \{\kap (\sig)\mid \sig <\del_{k-1}\}/D)\ .$$
Define $A^*=b_\mu [E^*]\cap \{\kap (\sig)\mid \sig
<\del_{k-1}\}$.  By Claim 1.25, then $A^*\in D$. From now we
continue as in 1.14 only using 1.10 in a fashion explained above
and at the final stage picking $A^{**}$ inside $D$.
\newline
\hfill $\square$

\medskip
\noindent
{\bf Remark 1.26.}\quad The use of Kinds$^*$ and not of Kinds in
1.21 (or actually in 1.22) is due only to our inability to extend
1.10 in order to include the case of a countable set.  Still in
view of 1.1 and also 1.23, 1.24, the first unclear case will not
be $\ome_1$ but rather $\ome_1+\ome_1$.

\section{Some Related Forcing Constructions}

In this section we like to show that (1) it is
impossible to remove SSH assumptions from Theorem
1.6; (2) the conclusion of Theorem 1.11 is optimal,
namely, starting with $\kap
=\bigcup_{n<\ome}\kap_n$,  $\kap_0
<\kap_1<\cdots <\kap_n <\cdots$  and $o(\kap_n)
=\kap^{+\del^n+1}_n +1$  we can construct a
model satisfying $2^\kap \ge\kap^{+\alp}$  for
every $\alp <\del^+$, where $\del$  as in 1.9
is a cardinal of uncountable cofinality; 
(3) the forcing construction for $\del$'s of
cofinality $\aleph_0$ will be given.  All these
results based on forcing of [Git1] and we sketch
them modulo this forcing.  

\medskip\noindent
{\bf Theorem 2.1}\quad {\sl Suppose that for
every $n<\ome$ $\{\alp <\kap\mid o(\alp)\ge\alp^{+n}\}$
is unbounded in $\kap$.  Then for every $\del
<\kap$  there is a cardinal preserving generic
extension such that it has at least $\del$ blocks
of principal indiscernibles $\langle\rho_{n,\nu}\mid
n<\ome,\nu <\del\rangle$ so that 
\begin{enumerate}
\item [(i)] $\rho_{n,\nu}<\rho_{n,\nu'}<\rho_{n+1,0}$
for every $n<\ome$, $\nu <\nu'<\del$
\item [(ii)]$\bigcup_{n<\ome}\rho_{n,\nu}=\kap$
for every $\nu <\del$\hb
and
\item [(iii)] $tcf\Big(\prodl_{n <\ome}
\rho^{+n+2}_{n,\nu}$, finite$\Big)=\kap^{++}$, for
every $\nu <\del$.
\end{enumerate}
}

\medskip\noindent
{\bf Proof.}\quad Without loss of generality we
can assume that $\del$  is a regular cardinal.
We pick an increasing sequence $\langle
\kap_n\mid n<\ome\rangle$ converging to $\kap$
so that for every $n<\ome$ $o(\kap_n)=\kap^{+n+2}_n+\del
+1$. Fix at each $n$  a coherent sequence of
extenders $\langle E^n_i\mid i\le\del\rangle$
with $E^n_i$  of the length $\kap^{+n+2}_n$.

We like to use the forcing of [Git1, Sec. 2] with
the extenders sequence $\langle E^n_\del \mid
n<\ome\rangle$  to blow power of $\kap$  to
$\kap^{++}$  together with extender based
Magidor forcing changing cofinality of the
principal indiscernible of $E^n_\del$ to $\del$ (for every
$n<\ome)$ simultaneously blowing its power to the
double plus.  We refer to M. Segal [Seg] or C.
Merimovich [Mer] for generalizations of the
Magidor forcing to the extender based Magidor
forcing. 

The definitions of both of these forcing notions
are rather lengthy and we would not reproduce
them here.  Instead let us emphasize what
happens with indiscernibles and why (iii) of the
conclusion of the theorem will hold.

Fix $n<\ome$.  A basic condition of [Git1, Sec. 2] is of
the form $\langle a_n,A_n,f_n\rangle$,
where $a_n$  is an order preserving function from
$\kap^{++}$  to $\kap_n^{+n+2}$  of cardinality
$<\kap_n$,  $A_n$  is a set of measure one for
the maximal measure of $rng a_n$ which is in
turn a measure of the extender $E^n_\del$  over
$\kap_n$.  The function of $f_n$  is an element
of the Cohen forcing over a $\kap^+$.  Each
$\alp\in dom a_n$  is intended to correspond to
indiscernible which would be introduced by the
measure $a_n(\alp)$ of $E^n_\del$.  In present
situation we force over the principal indiscernible
$\del_n$, i.e. one corresponding to the normal
measure of $E^n_\del$.  The extender based
Magidor forcing changes its cofinality to $\del$
and adds for every $\gam$ $\rho_n\le\gam\le\rho_n^{+n+2}$
a sequence $t_{n\gam}$  of order type $\del$
cofinal in $\rho_n$.  Actually, $t_{n,\rho_n^{+n+2}}(i)=
\rho_{n,i}^{+n+2}(i<\del)$, where $\langle\rho_{ni}|i<\del
\rangle$  is the sequence $t_{n\rho_n}$.  Now, if
$\gam <\rho_n^{+n+2}$  is produced by $a_n(\alp)$, then 
we connect $\alp$  with the sequence
$t_{n\gam}$ in addition to its connection 
with $\gam$.  Using standard arguments about
Prikry type forcing notions, it is not hard to
see that $cf\Big(\prodl_{n<\ome}\rho^{+n+2}_{n,i},{\rm
finite}\Big)=\kap^{++}$  for every
$i<\del$ as witnessed by $t_{n\gam}(i)'s$.\hfill
$\square$ 

\medskip\noindent
{\bf Remark 2.2}\quad Under the assumptions of
the theorem, one can obtain $2^\kap\ge\kap^{+\alp}$
for any countable $\alp$.  But we do not know
whether it is possible to reach uncountable
gaps.  See also the discussion in the final section. 

\medskip\noindent
{\bf Theorem 2.3}\quad {\sl Suppose that $\kap$ is a 
cardinal of cofinality $\ome,\del <\kap$  is a
cardinal of uncountable cofinality and for every
$n<\ome$  the set $\{\alp <\kap\mid
o(\alp)\ge\alp^{+\del^n}\}$ is unbounded in $\kap$.
{\it Then\/} for every $\alp <\del^+$  there is
cofinality preserving, not adding new bounded
subsets to $\kap$  extension satisfying $2^\kap\ge
\kap^{+\alp}$.
}

\medskip\noindent
{\bf Remark 2.4}  By the results of the previous
section, this is optimal if
$\alp\in [\bigcup_{n<\ome}\del^n,\del^+)$, at least if one forces
over the core model.

\medskip\noindent
{\bf Proof.}\quad Fix an increasing sequence
$\kap_0<\kap_1<\cdots <\kap_n<\cdots$
converging to $\kap$ so that each $\kap_n$
carries an extender $E_n$  of the length
$\kap_n^{+\del^n}$.  W.l. of $g$.  $\alp\ge
\bigcup_{n<\ome}\del^n$.  We use the Rado-Milner
Paradox (see K. Kunen [Kun, Ch. 1, ex. 20])  and
find $X_n\subset\alp (n\in\ome)$  such that
$\alp =\bigcup_{n<\ome} X_n$ and $otp
(X_n)\le\del^n$. W.l. of $g$. we can assume that
each $X_n$  is closed and $X_n\subseteq X_{n+1}(n<\ome)$.
Now the forcing similar to those of [Git1, 5.1]
will be applied.  Assign cardinals below $\kap$
to the cardinals $\{\kap^{+\bet +1}\mid 1\le\bet\le\alp\}$
as follows:  at level $n$  elements of the set
$\{\kap^{+\bet +1}\mid\bet +1\in X_n\}$  will correspond
to elements of the set $\{\kap^{+n+\gam+1}_n\mid
\gam <\del^n\}$.  

The next definition repeats 5.2 of [Git1] with
obvious changes taking in account the present
assignment.

\medskip\noindent
{\bf Definition 2.5}\quad The forcing noting
$\cal{P}(\alp)$  consists of all sequences
$\l\l A^{0\nu},A^{1\nu},F^\nu\r\mid\nu\le\alp\rangle$
so that 
\begin{enumerate}
\item [(1)] $\l\l\langle A^{0\nu},A^{1\nu}\rangle\mid\nu\le
\alp\rangle$  is as in 4.14 of [Git1].
\smallskip
\item [(2)] for every $\nu\le\alp$ $F^\nu$
consists of $p=\langle p_n\mid n<\ome\rangle$
and for every $n\ge\ell (p)$, $p_n=\langle
a_n,A_n,f_n\rangle$ as in 4.14 of [Git1] with the
following changes related only to $a_n$;
\end{enumerate}
\begin{enumerate}
\item [(i)] $a_n(\kap^{+\nu})= \kap_n^{+\varphi_n(\nu)}$
where $\varphi_n$  is some fixed in advance
order preserving function from successor
ordinals in $X_n$  to successor ordinals of  
$[n+2,\del^n)$.
\smallskip
\item [(ii)] only of cardinalities $\kap^{+\nu}$
for $\nu\in X_n\cap$  Successors can appear in
dom $a_n$.
\end{enumerate}

The rest of the argument repeats those of [Git1].

The following is a more general result that deals with all kinds
(i.e. elements of Kinds) of ordinals and not only with
$\del^n$'s.

\medskip
\noindent
{\bf Theorem 2.6.}\quad {\sl Let $\kap$  be a cardinal of
cofinality $\ome$  and $\del_0^{\ell_0}\cdots
\del_{k-1}^{\ell_{n-1}}\cdot\del\in\text{Kinds}\cap\kap$.
Suppose that for every $n<\ome$ the set
$\{\alp <\kap\mid o(\alp)\ge\alp^{+\del_0^{\ell_0}\cdots\del_{k-1}^{\ell_{k-1}}\cdot\del^n}\}$ is unbounded in $\kap$.
Then for every $\alp <\del_0^{\ell_0}\cdots\del_{k-1}^{\ell^{k-1}}
\cdot\del^+$  there is cofinality preserving, not adding new
bounded subsets to $\kap$  extension satisfying
$2^\kap\ge\kap^{+\alp}$.
}

Again, this is optimal by results of the
previous section, if
$$\alp\in [\bigcup_{n<\ome}\del_0^{\ell_0}\cdots
\del_{k-1}^{\ell_{k-1}}\cdot\del^n\ ,\ \del_0^{\ell_0}
\cdots \del_{k-1}^{\ell_{k-1}}\cdot\del^+)$$
at least if one forces over the core model in case $\del
=\aleph_1$. The construction is parallel to those of 2.3, only we
use the following version of Rado-Milner Paradox:

For every $\alp\in [\bigcup_{n<\ome}\del_0^{\ell_0}\cdots
\del_{k-1}^{\ell_{k-1}}\cdot\del^n$, $\del_0^{\ell_0}\cdots
\del_{k-1}^{\ell_{k-1}}\cdot\del^+)$  there are $X_n\subseteq
\alp (n<\ome)$ such that $\alp =\bigcup_{n<\ome}X_n$ and
$otp (X_n)\le\del_0^{\ell_0}\cdots\del_{k-1}^{\ell_{k-1}}
\cdot\del^n$.
\newline
\hfill $\square$

Under the same lines we can deal with gaps of
size of a cardinal of countable cofinality below
$\kap$.  Thus the following result which together
with the results of the previous section provides
the equiconsistency holds:

\medskip\noindent
{\bf Theorem 2.7}\quad {\sl Suppose that $\kap$
is a cardinal of cofinality $\ome$  and $\del
<\kap$ is a cardinal of cofinality $\ome$  as
well.  Assume that for every $\tau <\del$  the
set $\{\alp <\kap\mid o(\alp)\ge\alp^{+\tau}\}$
is unbounded in $\kap$.  {\it Then\/} for every
$\alp <\del^+$  there are cofinalities preserving,
not adding new bounded subsets to $\kap$
extension satisfying $2^\kap\ge\kap^{+\alp}$. 
}

The proof is similar to those of 2.3.  Only notice that we can
present $\alp$ as an increasing union of sets $X_n(n<\ome)$  with
$|X_n|<\del$  since $\alp <\del^+$, $cf\del =\ome$ and
there is a function from $\del$  onto $\alp$.

\section{Concluding Remarks and Open Questions} 

Let us first summarize in the table below the situation under
$SSH_{<\kap}$  (i.e.  for every singular $\mu <\kap$  $pp\mu=\mu^+$)
assuming that $2^\kap\ge\kap^{+\del}$  for some $\del$, where
$\kap$  as usual here in a strong limit cardinal of cofinality
$\aleph_0$.  For $\del =\aleph^\ell_1$, for $2\le\ell <\ome$,
in the cases dealing with ordinals in ${\rm Kinds}\bks{\rm Kinds}^*$ 
we assume in addition that there is no
measurable of the core model between $\kap$  and
$\kap^{+\del}$.

\begin{center}
{\scriptsize
\begin{tabular}{|c|c|c|c|}
\hline
$\delta = 2$&\multicolumn{3}{c|}{$o(\kappa)
= \kappa^{++}$}\\
&\multicolumn{3}{c|}{or}\\
&\multicolumn{3}{c|}{$\forall n<\ome\{\alpha <
\kappa | o(\alpha) \ge \alpha^{+n}\}$}\\
&\multicolumn{3}{c|}{is unbounded in
$\kappa$}\\
\hline
$2 < \delta < \aleph_0$&\multicolumn{3}{c|}
{$o(\kappa) = \kappa^{+\delta} + 1$}\\
&\multicolumn{3}{c|}{or}\\
&\multicolumn{3}{c|}
{$\forall n<\ome\{\alpha <\kappa |
o(\alpha)\ge \alpha^{+n}\}$}\\
&\multicolumn{3}{c|} {is unbounded in
$\kappa$}\\
\hline
&$cf|\delta| =\aleph_0$&\multicolumn{2}{c|} {$\forall \tau
< |\delta| \ \{\alpha < \kappa | o(\alpha)
\ge \alpha^{+\tau}\}$}\\
&&\multicolumn{2}{c|}{ is unbounded in
$\kappa$}\\
\cline{2-4}
&& $\delta$ is a cardinal&
$o(\kappa) \ge \kappa^{+\delta + 1} + 1$\\
$\kappa > \delta \ge \aleph_0$&$cf |\delta |
> \aleph_0$
&&or\\
&&&$\{\alpha < \kappa | o(\alpha) \ge
\alpha^{+\delta + 1} + 1\}$\\
&&&is unbounded in $\kappa$\\
\cline{3-4}
&&$\delta = |\delta|^\ell$,&$o(\kappa) \ge
\kappa^{+|\delta|^\ell + 1} + 1$\\
&&for some&or\\
&&$1 < \ell < \omega$&$\{\alpha < \kappa |
o(\alpha) \ge^{+|\delta|^\ell +
1} + 1\}$\\
&&&is unbounded in $\kappa$\\
\cline{3-4}
&&$\delta\ge\bigcup\limits_{\ell < \omega}
|\delta|^\ell$&$\forall \ell < \omega
\{\alpha < \kappa|o(\alpha)\ge\alpha^{+|\delta|^\ell}\}$\\
&&&{is unbounded in $\kappa$}\\
\cline{3-4}
&&$\delta_0^{\ell_0}\cdots\delta_{k-1}^{\ell_{k-1}}\cdot
\delta^\omega_k{\le}\del{<}\del_0^{\ell_0}\cdots
\delta_{k-1}^{\ell_{k-1}}\cdot\delta^+_k$&$\forall
n{<}\ome\{\alp{<}\kappa\mid o(\alp){\ge}
\alp^{+\delta_0^{\ell_0}\cdots\delta_{k-1}^{\ell_{k-1}}
\cdot\delta^n_k}\}$\\
&&for some $\delta_0^{\ell_0}\cdots
\delta_{k-1}^{\ell_{k-1}}\cdot\del_k\in$ Kinds&is
unbounded in $\kappa$\\
\cline{3-4}
&&$\delta_0^{\ell_0}\cdots\delta_k^{\ell_k}\le\del
<\del_0^{\ell_0}\cdots\del_k^{\ell_k}
\cdot\ome_1$&$o(\kap)\ge\kap^{+\del_0^{\ell_0}\cdots
\del_k^{\ell_k}+1}+1$\\
&&for some $\del_0^{\ell_0}\cdots\del_k^{\ell_k}\in$
Kinds&or\\
&&&$\{\alp <\kap\mid
o(\alp)\ge\alp^{+\del_0^{\ell_0}\cdots\del_k^{\ell_k}+1}
+1\}$\\
&&&is unbounded in $\kap$\\
\hline
$\del\ge\kap$&\multicolumn{3}{c|}
{$\forall\tau <\kap \{\alp
<\kap\mid o(\alp)\ge \alp^{+\tau}\}$}\\
&\multicolumn{3}{c} {is unbounded in $\kappa$}\\
\hline
\end{tabular}
}
\end{center}

The proofs are spread through the papers [Git1,2,3,4,5],
[Git-Mag], [Git-Mit] and the present paper.
The forcing constructions in these papers give
GCH below $\kap$.  	 

Let us finish with some open problems.  

\medskip\noindent
{\bf Question 1.}\quad Let $a$ be a countable
set of regular cardinals.  Does ``$|pcf
a|>|a|=\aleph_0$" imply an inner model with a
strong cardinal? 

In view of 1.10, it is natural to understand the
situation for countable $a$.  Recall that the
consistency of ``$|pcfa|>|a|$" is unknown and it
is a major question of the cardinal arithmetic. 

The next question is more technical.

\medskip\noindent
{\bf Question 2.}\quad Can the assumption that
there are no measurables in the core model
between $\kap$  and $2^\kap$ be removed in 1.11?

It looks like this limitation is due only to the
weakness of the proof. But probably there is a
connection with ``$|pcf a|>|a|"$. The simplest
unclear case is $2^\kap\ge\kap^{+\ome^2_1}$. 

The situation without $SSH_{<\kap}$  is unclear.
In view of 2.1 probably weaker assumptions then
those used in the case of $SSH_{<\kap}$  may
work.  A simplest question in this direction is
as follows.

\medskip\noindent
{\bf Question 3.}\quad Is ``$\{\alp\mid
o(\alp)\ge\alp^{+n}\}$ unbounded in $\kap$  for
each $n<\ome$" sufficient for ``$\kap$  strong
limit, $cf\kap =\aleph_0$  and
$2^\kap\ge\kap^{+\ome_1}$"?

If the answer is affirmative, then the
construction will require a new forcing with
short extenders, which will be interesting by
itself.   We then conjecture
that the same assumption will work for arbitrary
gap as well.

For uncountable cofinalities (i.e. $cf\kap
>\aleph_0$), as far as we are concerned with
consistency strength, the only unknown case is
the case of cofinality $\aleph_1$. We restate a
question of [Git-Mit]:

\medskip\noindent
{\bf Question 4.}\quad What is the exact
strength of ``$\kap$  is a strong limit, $cf\kap
=\aleph_1$  and $2^\kap\ge\lam$  for a regular
$\lam>\kap^+$?

It is known that the strength lies between $o(\kap)=\lam$
and $o(\kap)=\lam +\ome_1$, see [Git-Mit].


\begin{thebibliography}{2}
\bibitem[Git5]{Git5} M Gitik, There is no bound
for the power of the first fixed point.
\end{thebibliography}
\end{document}